\documentclass[12pt]{amsart}

\usepackage{color}

 \definecolor{formula}{rgb}{0.2,0.5,0.3}
 
\usepackage{tikz}
\usepackage{verbatim}

\newcommand{\elle}[1]{L^{#1}(\Omega)}
\newcommand{\w }[1]{W^{1,{#1}}_{0}(\Omega)}
\newcommand{\N}{{I\!\!N}}
\newcommand{\R}{{I\!\!R}}
\newcommand{\un}{u_{n}}
 
\newcommand{\norma}[2]{\|#1\|_{\lower 4pt \hbox{$\scriptstyle #2$}}}

\newcommand{\io}{\int\limits_{\Omega}}
\newcommand{\ik}{\int\limits_{\{k\leq|u_n|\}}}
\newcommand{\iu}{\int\limits_{\{1\leq|u_n|\}}}
 \newcommand{ \dive}{{\rm div}}
 
  \newcommand{ \weak}{\rightharpoonup}

\long\def\salta#1{\relax}
\newcommand{\arrstre}{\renewcommand{\arraystretch}{2}}

\newcommand{\rife}[1]{(\ref{#1})}
\newtheorem{example}{Example}[section]

\newtheorem{theo}[example]{\sc Theorem}

\def\picture #1 by #2 (#3){
  \vbox to #2{
    \hrule width #1 height 0pt depth 0pt
    \vfill
    \special{picture #3} 
    }
  }

\def\scaledpicture #1 by #2 (#3 scaled #4){{
  \dimen0=#1 \dimen1=#2
  \divide\dimen0 by 1000 \multiply\dimen0 by #4
  \divide\dimen1 by 1000 \multiply\dimen1 by #4
  \picture \dimen0 by \dimen1 (#3 scaled #4)}
  }

\begin{document}

 
\title[$W_0^{1,1}$ solutions in a problem with degenerate coercvity]{$W_0^{1,1}$ solutions in some borderline cases of  elliptic equations  with degenerate coercivity }

\author{Lucio Boccardo, Gisella Croce}
\address{L.B. -- Dipartimento di Matematica, ``Sapienza'' Universit\`{a} di Roma,
P.le A. Moro 2, 00185 Roma (ITALY)}
\email{boccardo@mat.uniroma1.it}
\address{G.C. -- Laboratoire de Math\'ematiques Appliqu\'ees du Havre, Universit\'e du Havre,
25, rue Philippe Lebon, 76063 Le Havre (FRANCE)}
\email{gisella.croce@univ-lehavre.fr}



\begin{abstract}
We study a  degenerate elliptic equation, proving existence results of   distributional solutions in some borderline cases.
\end{abstract}

\maketitle

\rightline{{\sl Bernardo, come vide li occhi miei ...}}\rightline{(Dante, Paradiso XXXI)}

\section*{Acknowledgements.} 
This paper contains 
developments of  the results presented  by the first author   at IX WNDE (Jo\~ao Pessoa, 18.9.2012), 
  {\sl on the occasion of the sixtieth  birthday of Bernhard Ruf\,}'',
and it is dedicated to him.

\section{Introduction}

 {
 The Sobolev space $\w2$ is the natural
functional framework (see \cite{br-fourier},  \cite{ll}) to find weak solutions
of nonlinear elliptic problems
 of the following type
\begin{equation}
   \left\{ \arrstre
      \begin{array}{cl}
  \displaystyle      -\dive \left(\frac{a(x)\nabla   u}{(1+|u|)^\theta}\right) = f, & \mbox{in $\Omega$;}\\
         u=0 , & \mbox{on $\partial\Omega$,}
     \end{array}
   \right.
\label{p}
\end{equation}
where  
the
function $f$ belongs to the dual space of $\w2$,
$\Omega$ is a bounded, open subset of $\R^N$, with $N > 2$,
 $\theta$ is a real number such that
\begin{equation}
     0 \leq \theta\leq 1\,,
\label{t1}
\end{equation}
and $a:\Omega \to
\R$ is a  measurable function  satisfying the following conditions:
\begin{equation}
  {\alpha}\leq  a(x) \leq \beta\,,
\label{a1}
\end{equation}
for almost every $x \in \Omega$, where $\alpha$ and $\beta$ are positive constants.  
{The main difficulty to use the general results of \cite{br-fourier},  \cite{ll}
is the fact that}
$$
A(v) = \displaystyle -\dive \left(\frac{a(x)\nabla   v}{(1+|v|)^\theta} \right), 
$$
is not coercive.
{Papers \cite{bdo}, \cite{bb} and \cite{b} deal with the existence and summability of solutions to  problem
(\ref{p}) 
if $f\in \elle m$, for some $m\geq 1$}.
 
Despite the lack of coercivity of the
differential operator 
$A(v)$
appearing in  problem \rife{p}, in the papers {\cite{bdo}, \cite{bb} and \cite{ABFOT}}, the authors prove the following existence results of solutions of problem \rife{p}, under assumption \rife{a1}:
}  
\begin{itemize}
\item[A)] a weak solution  $u\in\w 2
\cap \elle\infty$, if $m > \frac N2$ and   \rife{t1} holds true;
\item[B)]
a weak solution 
$u\in  \w 2\cap\elle  {m^{**}(1-\theta)}$, where $m^{**}=(m^{*})^{*}=\frac {mN}{N-2m}$, if  
$$
0<\theta <1,\quad
{\frac{2N}{N+2-\theta(N-2)} }
 \leq m < \frac N2\,;
$$
\item[C)]
a distributional solution  $u$ in $\w q$, 
$
\displaystyle q = \frac{Nm(1-\theta)}{N-m(1+\theta)} < 2\,,
$ 
 if  
$$
{\frac1{N-1} \leq} \theta <1,\quad
{
\frac{N}{N+1-\theta(N-1)} 
}
< m < \frac{2N}{N+2-\theta(N-2)}\,.
$$
\item[D)]
an entropy solution $u \in M^{m^{**}(1-\theta)}$, with $|\nabla u| \in M^{q}(\Omega)$,
for $1\leq m\leq\max \left\{1,\frac{N}{N+1-\theta(N-1)}\right\}$.
\end{itemize}
 The borderline case $\theta=1$ was studied in \cite{b}, proving the existence of a solution $u \in W^{1,2}_0(\Omega)\cap L^p(\Omega)$ for every $p<\infty$. The case where the source is $\frac{A}{|x|^2}$ was analyzed too. 
  
About the different notions of solutions mentioned above, we recall that the notion of entropy solution was introduced in \cite{BBGGPV}.
Let 
\begin{equation}
T_k(s)= \left\{
\begin{array}{ll}
s
& \mbox{if $|s|\leq k$,}
\\
k\frac{s}{|s|}
& \mbox{if $|s|> k$.} 
\end{array}
\right.
\label{trunc}
\end{equation}  
Then $u$ is an entropy solution to problem (\ref{p})  if $T_k(u) \in {\w2}$ for every $k>0$
and
$$
\io \frac{a(x)\nabla u}{(1+|u|)^{\theta}}\cdot \nabla T_k(u-\varphi)\leq \io f\,T_k(u-\varphi)\,,\quad\quad
\forall\,\varphi \in {\w2}\cap L^{\infty}(\Omega)\,.
$$
Moreover, we say that 
$u$ is a distributional solution of \rife{p}  if
\begin{equation}
\io {\frac{a(x)\,\nabla   u}{(1+|u|)^\theta} \cdot \nabla  \varphi} = \io {f\,\varphi}
\,,
\qquad
\forall\; \varphi \in C^\infty_0(\Omega)\,.
\label{distrib}
\end{equation} 
The figure below can help to summarize the previous results, where the name of a given region corresponds to the results that we have just cited. 
\begin{center}
\includegraphics[width=6cm]{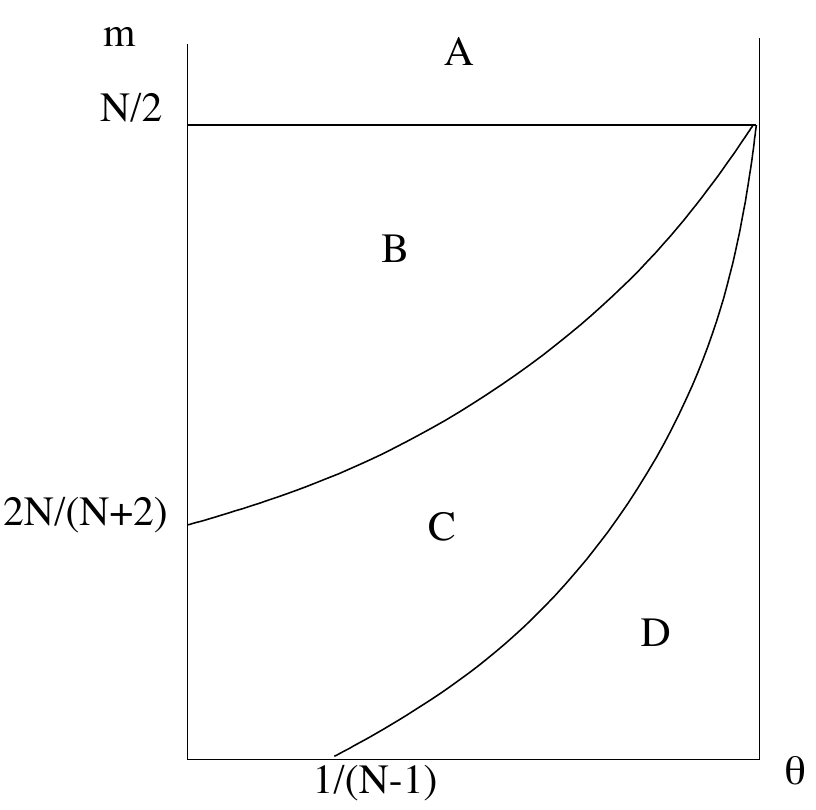}  
\end{center}

Our results are the following.
 
 \begin{theo}\label{th1}\sl
Let $f$ be a function in $\elle m$. 
 Assume \rife{a1} and
 \begin{equation}
\label{m-teta}
m=\frac{N}{N+1-\theta(N-1)},
\quad \frac{1}{N-1}<\theta<1
\end{equation}
Then there exists  a $\w 1$ distributional solution to problem  \rife{p}.
\end{theo}
Observe that this case corresponds to the curve between the regions $C$ and $D$ of the figure.
Note that $m>1$ if and only if $\theta>\frac1{N-1}$. 

In the following result too, we will prove the existence of a $\w1$ solution.
\begin{theo}\label{th4}
Let $f$ be a function in $\elle {m}$. 
 Assume \rife{a1},
 $f\log(1+|f|)\in\elle1$ and
$\theta=\frac{1}{N-1}$.
Then there exists a $\w1$ distributional solution of \rife{p}.
 \end{theo}

We end our introduction just mentioning a uniqueness result of solutions to problem \rife{p} can be found in
\cite{Po}. 

Moreover, in \cite{bb, BCOManuscripta, croce_rendiconti, BCOLincei}
it was showed that the presence of a lower order term has a regularizing effect on the existence and regularity of the solutions.

To prove our results, we will work by approximation, using the following sequance of problems:
\begin{equation}
   \left\{ \arrstre
      \begin{array}{cl}
         \displaystyle -{\rm div}  \left(\frac{a(x)\nabla   u_n}{(1+|u_n|)^{\theta}}\right) = T_n(f), & \mbox{in $\Omega$;}\\
         u_n = 0, & \mbox{on $\partial\Omega$.}
   \end{array}
   \right.
\label{pn}
\end{equation} 
The existence of weak solutions $u_n \in \w 2 \cap L^{\infty}(\Omega)$  to problem
\rife{pn}
is due to \cite{bdo}.

\section{$W^{1,1}_0(\Omega)$ solutions}

In the sequel $C$ will denote a constant depending on $\alpha, N, {\rm meas}(\Omega), \theta$ and the $L^m(\Omega)$ norm of the source $f$.

We are going to prove Theorem \ref{th1}, that is, the existence of a solution to problem \rife{p}
in the case where
$m=\frac{N}{N+1-\theta(N-1)}$ and
 $\frac{1}{N-1}<\theta<1$.
Note that $m>1$ if and only if $\theta>\frac1{N-1}$. 
\begin{proof}(of Theorem \ref{th1})
We consider $T_k(u_n)$ as a test function in (\ref{pn}): then
\begin{equation}
\label{tk1}
\io|\nabla T_k(u_n)|^2\leq
\frac{(1+k)^\theta\,\norma{f}{\elle1}}{\alpha}
\end{equation}
by assumption (\ref{a1}) on $a$.

Choosing  
$
 [(1+|u_n|)^p - 1] {\rm sign}(u_n)
$, for $p=\theta - \frac{1}{N-1}$,
as a test function in \rife{pn}
we have, by H\"older's inequality on the right hand side and assumption (\ref{a1}) on the left one
\begin{equation}\label{inizio}
\alpha\,p\io\left(\frac{|\nabla u_n|}{(1+|u_n|)^\frac{N}{2(N-1) }}\right)^2
\leq \io|f|[(1+|u_n|)^p - 1]  
\leq\norma{f}{\elle m}\bigg[\io[(1+|u_n|)^p - 1] ^{m'} \bigg]^\frac1{m'}.
\end{equation}
The   Sobolev embedding used on the left hand side   implies
$$
\bigg[
\io\bigg\{(1+|u_n|)^\frac{N-2}{2(N-1) }-1\bigg\}^\frac{2N}{N-2}
\bigg]^\frac2{2^*}
\leq
C\bigg[\io[(1+|u_n|)^p - 1] ^{ m'} \bigg]^\frac1{m'}.
$$
We observe that 
$
\frac{N-2}{2(N-1) }\frac{2N}{N-2}=pm'
$; moreover
$
\frac2{2^*}>\frac1{m'}$, { since } $m<\frac N2$.
Therefore the above inequality implies that  
\begin{equation}\label{r}
\io|u_n|^\frac{N}{ N-1}\leq C.
\end{equation} 
One deduces that
\begin{equation}\label{l}
 \io\frac{|\nabla u_n|^2}{(1+|u_n|)^\frac{N}{N-1}}
\leq  C
\end{equation}
from (\ref{r}) and (\ref{inizio}).
Let 
$v_n=\frac{2(N-1)}{N-2}(1+|u_n|)^{\frac{N-2}{2(N-1)}}{\rm {sign}} (u_n)$. Estimate \rife{l} is equivalent to say that $\{v_n\}$ is a bounded sequence in
${\w2}$; therefore, up to a subsequence, there exists $v \in {\w2}$
such that $v_n\weak v$ weakly in ${\w2}$ and a.e. in $\Omega$. 
If we define the function
 $u=\Big(\left[\frac{2(N-1)}{N-2}|v|\right]^{\frac{2(N-1)}{N-2}}-1\Big){\rm{sign}}(v)$, the weak convergence of $\nabla v_n\weak \nabla v$
means that
\begin{equation}
\frac{\nabla u_n}{(1+|u_n|)^{\frac{N}{2(N-1)}}}
\weak \frac{\nabla u}{(1+|u|)^{\frac{N}{2(N-1)}}}\quad \textnormal{weakly in}\,\,L^2(\Omega)\,.
\label{convergenza_v_n}
\end{equation}
Moreover, the Sobolev embedding for $v_n$ implies that $u_n\to u$ in $L^s(\Omega)$, for every $1\leq s<{\frac{N}{N-1}}$.

Hold\"er's inequality with exponent 2
applied to
$$
\io|\nabla u_n|=
\io
\frac{|\nabla u_n|}{(1+|u_n|)^\frac{N}{2N-2 }}
(1+|u_n|)^\frac{N}{2N-2}
$$
gives 
\begin{equation}\label{1}
\io|\nabla u_n|\leq 
 C,
\end{equation}
due to 
 (\ref{r}) and (\ref{l}). 
We are now going to estimate $\displaystyle \ik |\nabla u_n|$. By using 
$
 [(1+|u_n|)^p - (1+k)^p]^+ {\rm sign}(u_n)
$
as a test function  in \rife{pn},
we have 
$$
\ik\frac{|\nabla u_n|^2}{(1+|u_n|)^\frac{N}{N-1}}
\leq C\bigg[\ik|f|^m\bigg]^\frac1m
\bigg[\io(1+|u_n|)^\frac N{N-1} \bigg]^\frac1{m'},
$$
which implies, by (\ref{r}),
\begin{equation}\label{iniziok}
\ik\frac{|\nabla u_n|^2}{(1+|u_n|)^\frac{N}{N-1}}
\leq\,C\bigg[\ik|f|^m\bigg]^\frac1m\,.
\end{equation}
H\"older's inequality, estimates (\ref{r}) and (\ref{iniziok})
on 
$$
\ik|\nabla u_n|=
\ik
\frac{|\nabla u_n|}{(1+|u_n|)^\frac{N}{2N-2 }}
(1+|u_n|)^\frac{N}{2N-2 }
\leq\,C\bigg[\ik|f|^m\bigg]^\frac1{2m}\,,
$$
give
\begin{equation}\label{1k}
\ik|\nabla u_n|=
\ik
\frac{|\nabla u_n|}{(1+|u_n|)^\frac{N}{2N-2 }}
(1+|u_n|)^\frac{N}{2N-2 }
\leq\,C\bigg[\ik|f|^m\bigg]^\frac1{2m}\,.
\end{equation}
Thus, for every   measurable subset $E$,   due to \rife{tk1} and \rife{1k},
we have
\begin{equation}
\label{malaga}
\begin{array}{l}
\displaystyle
\int\limits_E\bigg|\frac{\partial u_n}{\partial x_{i}}\bigg|
\leq
\int\limits_E|\nabla u_n|
\leq
\int\limits_E|\nabla  T_{k}(u_n)|
+
\ik|\nabla u_n|\, 
\\
\displaystyle
 \leq{\rm meas}(E)^\frac1{2}
\left[\frac{(1+k)^\theta\,\norma{f}{\elle1}}{\alpha} \right]^\frac12
+
C\bigg[\ik|f|^m\bigg]^\frac1{2m}\,.
 \end{array}
\end{equation}
Now we are going to prove that 
$u_n$ weakly converges to $u$ in $\w1$ 
following \cite{BCOManuscripta}.
Estimates \rife{malaga} and \rife{r} imply that the sequence $\{\frac{\partial{u_n}}{\partial x_{i}}\}$ is equiintegrable.
By Dunford-Pettis theorem, and up to subsequences, there exists $Y_{i}$ in $\elle1$ such that $\frac{\partial{u_n}}{\partial x_{i}}$ weakly converges to $Y_{i}$ in $\elle1$. Since $\frac{\partial{u_n}}{\partial x_{i}}$ is the distributional partial derivative of ${u_n}$, we have, for every $n$ in $\N$,
$$
\io \frac{\partial{u_n}}{\partial x_{i}}\,\varphi = -\io {u_n}\,\frac{\partial\varphi}{\partial x_{i}}\,,
\quad
\forall\; \varphi \in C^{\infty}_{0}(\Omega)\,.
$$
We now pass to the limit in the above identities, using that $\partial_{i}{u_n}$ weakly converges to $Y_{i}$ in $\elle1$, and that ${u_n}$ strongly converges to $u$ in $\elle1$: we obtain
$$
\io Y_{i}\,\varphi = -\io u\,\frac{\partial\varphi}{\partial x_{i}}\,,
\quad
\forall\; \varphi \in C^{\infty}_{0}(\Omega)\,.
$$
This implies that $Y_{i} = \frac{\partial u}{\partial x_{i}}$, and this result is true for every $i$. Since $Y_{i}$ belongs to $\elle1$ for every $i$, $u$ belongs to $\w1$.

We are now going to pass to the limit in problems (\ref{pn}). For the limit of the left hand side, 
it is sufficient to observe that $\frac{\nabla u_n}{(1+|u_n|)^{\frac{N}{2(N-1)}}}\weak\frac{\nabla u}{(1+|u|)^{\frac{N}{2(N-1)}}}$
weakly in $L^2(\Omega)$ due to (\ref{convergenza_v_n}) and that $\left|a(x)\nabla \varphi\right|$ is bounded.
\end{proof}

We  prove Theorem \ref{th4},
that is, the existence of a $W^{1,1}_0(\Omega)$ solution in the case where
$\theta=\frac{1}{N-1}$ and   $f\log(1+|f|)\in\elle1$.

\begin{proof}(of Theorem \ref{th4})
Let $k\geq 0$ and take $ [\log(1+|\un|)-\log(1+k)]^+{\rm sign}(\un)$, as a test function in problems (\ref{pn}).
By assumption (\ref{a1}) on $a$ one has
$$\alpha
\ik\frac{|\nabla\un|^2\; }
{(1+|\un|)^{\theta+1}}
\leq
\ik|f|\log(1+|\un|)\,.
$$ 
We now use the following inequality on the left hand side:
 \begin{equation}
\label{bar}
a\log(1+b)
\leq \frac a\rho\log\left(1+\frac a\rho\right)+(1+b)^\rho
\end{equation}
where $a,b$ are positive real numbers and $0<\rho<\frac{N-2}{N-1}$.  This gives, for any $k\geq 0$
\begin{equation}
\label{camino0}
\alpha
\ik\frac{|\nabla\un|^2\;  }
{(1+|\un|)^{\theta+1}  }
\leq\ik\frac{|f|}{\rho}\log\bigg(1+\frac{|f|}{\rho}\bigg)
+\ik(1+|\un|)^\rho\,.
\end{equation}
In particular, for $k\geq 1$ we have
\begin{equation}
\label{camino}
\frac{\alpha}{2^{\theta+1}}
\ik\frac{|\nabla\un|^2\;  }
{|\un|^{\theta+1}  }
\leq\ik\frac{|f|}{\rho}\log\bigg(1+\frac{|f|}{\rho}\bigg)
+2^\rho\ik|\un|^\rho\,.
\end{equation}
Writing the above inequality for $k=1$ and using the Sobolev inequality on the left hand side, one has
$$
\bigg[\io
(|\un|^\frac{1-\theta}{2}-1)_+^{2^*} 
\bigg]^\frac2{2^*}
\leq C\iu\frac{|f|}{\rho}\log\bigg(1+\frac{|f|}{\rho}\bigg)
+C\iu|\un|^\rho\,,
$$
which implies that     
$$
\bigg[\io
 |\un|^\frac{(1-\theta)2^*}{2} 
\bigg]^\frac1{2^*}
\leq
C+C
\sqrt{\;\iu\frac{|f|}{\rho}\log\bigg(1+\frac{|f|}{\rho}\bigg)\;}
+C
\sqrt{\;\io|\un|^\rho\;}\,.
$$
By using H\"older's inequality with exponent $\frac{(1-\theta)2^*}{2\rho}$ on the last term of the right hand side, we get 
$$
\bigg[\io
 |\un|^\frac{(1-\theta)2^*}{2} 
\bigg]^\frac1{2^*}
\leq
C+C
\sqrt{\;\iu\frac{|f|}{\rho}\log\bigg(1+\frac{|f|}{\rho}\bigg)\;}
+
C
\bigg[\io|\un|^\frac{(1-\theta)2^*}{2}
\bigg]^\frac{\rho}{(1-\theta)2^*}\,.
$$
By the choice of $\rho$, this inequality implies that
\begin{equation}\label{stima_norma_lebesgue_u_n} 
\io|\un|^\frac{(1-\theta)2^*}{2}\leq\,C\,.
\end{equation}
Inequalities (\ref{stima_norma_lebesgue_u_n}) and (\ref{camino0}) written for $k=0$
imply that
$\{v_n\}=\{\frac{2}{1-\theta}(1+|u_n|)^{\frac{1-\theta}{2}}{\rm {sign}}(u_n)\}$ is a bounded sequence in
${\w2}$, as in the proof of Theorem \ref{th1}. Therefore, up to a subsequence there exists $v \in {\w2}$
such that $v_n\weak v$ weakly in ${\w2}$ and a.e. in $\Omega$. 
 Let $u=\left\{\left[\frac{1-\theta}{2} v\right]^{\frac{2}{1-\theta}}-1\right\}{\rm{sign}}(v)$;  the weak convergence of $\nabla v_n\weak \nabla v$
means that
\begin{equation}
\frac{\nabla u_n}{(1+|u_n|)^{\frac{\theta+1}{2}}}
\weak \frac{\nabla u}{(1+|u|)^{\frac{\theta +1}{2}}}\quad \textnormal{weakly in}\,\,L^2(\Omega)\,.
\label{convergenza_v_n_4}
\end{equation}
Moreover, the Sobolev embedding for $v_n$ implies that $u_n\to u$ in $L^s(\Omega), s<{\frac{N}{N-1}}$.

By  (\ref{tk1})
one has
$$
\io|\nabla\un|=
\io|\nabla T_1(\un)|+\iu|\nabla\un|
\leq
C+
\iu
\frac{|\nabla\un|}{|\un|^\frac{\theta+1}{2}}
|\un|^\frac{\theta+1}{2}\,.
$$
H\"older's inequality on the right hand side, and estimates
 (\ref{camino}) written with $k=1$ and 
(\ref{stima_norma_lebesgue_u_n})
imply that the sequence $\{\un\}$ is bounded in $\w1$.

Moreover, due to \rife{camino}
$$
\ik|\nabla\un|
=\ik
\frac{|\nabla\un|}{|\un|^\frac{\theta+1}{2}}
|\un|^\frac{\theta+1}{2}
\leq
$$
$$
\leq {C}
\sqrt{
\ik\frac{|f|}{\rho}\log\bigg(1+\frac{|f|}{\rho}\bigg)
+\ik|\un|^\rho
}\,.
$$
For every   measurable subset $E$,   the previous inequality and (\ref{tk1}) imply
$$\int\limits_E\bigg|\frac{\partial u_n}{\partial x_{i}}\bigg|
\leq
\int\limits_E|\nabla u_n|
\leq
\int\limits_E|\nabla  T_{k}(u_n)|
+
\ik|\nabla u_n|\, 
$$
$$
 \leq C\,{\rm meas}(E)^{\frac1{2}}
(1+k)^{\frac{\theta}{2}}
+
C\,\sqrt{
\ik\frac{|f|}{\rho}\log\bigg(1+\frac{|f|}{\rho}\bigg)
+\ik|\un|^\rho
}\,.
$$
Since $\rho<\frac{(1-\theta)2^*}{2}$, by using H\"older's inequality on the last term and
estimate (\ref{stima_norma_lebesgue_u_n}),
one has
$$
\int\limits_E\bigg|\frac{\partial u_n}{\partial x_{i}}\bigg|
\leq C\,{\rm meas}(E)^{\frac1{2}}
(1+k)^{\frac{\theta}{2}}+
$$
$$+
C\,\sqrt{
\ik\frac{|f|}{\rho}\log\bigg(1+\frac{|f|}{\rho}\bigg)
+{\rm meas}(\{|u_n|\geq k\})^{1-\frac{2\rho}{(1-\theta)2^*}}
}\,. 
$$
One can argue as in the proof of Theorem \ref{th1} to deduce that $u_n\to u$
weakly in $W^{1,1}_0(\Omega)$.

To pass to the limit in problems (\ref{pn}),  as in the proof of Theorem \ref{th1}, it is sufficient to observe that $\frac{\nabla u_n}{(1+|u_n|)^{\frac{N}{2(N-1)}}}\weak\frac{\nabla u}{(1+|u|)^{\frac{N}{2(N-1)}}}$
weakly in $L^2(\Omega)$, due to (\ref{convergenza_v_n_4}). 
\end{proof}


\begin{thebibliography}{99}
 
 
 
 \bibitem{ABFOT} A. Alvino, L. Boccardo, V. Ferone, L. Orsina, G. Trombetti:
 {  Existence results for nonlinear elliptic equations 
  with degenerate coercivity};


    Ann. Mat. Pura Appl.  182 (2003),  53--79. 
    
\bibitem{BBGGPV} P. Benilan, L. Boccardo, T. Gallou\"et, R. Gariepy, M. Pierre, and J.L. Vazquez:
An $L^1$ theory of existence and uniqueness of nonlinear elliptic equations; {Ann. Scuola Norm.
Sup. Pisa Cl. Sci.} {22} (1995),  240--273.

 \bibitem{b}
L. Boccardo:  Some elliptic problems with  
degenerate coercivity; Adv. Nonlinear Stud. 6 (2006),   1--12.

  \bibitem{bb}
L. Boccardo, H. Brezis:   Some remarks on a class of
elliptic equations with degenerate coercivity; {  Boll.
Unione Mat. Ital.}  {6} (2003),  521--530.


\bibitem{BCOManuscripta}
L. Boccardo, G. Croce, L. Orsina:  {Nonlinear degenerate elliptic
problems with $W^{1,1}_0(\Omega)$ solutions}; Manuscripta Mathematica 137 (2012), 419--439.


\bibitem{BCOLincei}
L. Boccardo, G. Croce, L. Orsina:  {  A semilinear problem with a
$W^{1,1}_0(\Omega)$ solution};   Rend. Lincei. Mat. Appl. 23 (2012), 97-103.

\bibitem{bdo} L. Boccardo, A. Dall'Aglio, L. Orsina:
  {   Existence and regularity results 
for some  elliptic equations with degenerate coercivity }; Atti Sem. Mat. Fis. Univ.
Modena 46 (1998),   51--81.

\bibitem{BGa} L. Boccardo, T. Gallou\"et: {  Nonlinear elliptic equations with
right-hand side measures}; {Comm. Partial Differential Equations} 17 (1992), 641--655.


\bibitem{BGa2012} L. Boccardo, T. Gallouet:
$W_0^{1,1}$ solutions in some borderline cases of
 Calderon-Zygmund theory; J. Differential Equations 253 (2012), 2698--2714.
 
 \bibitem{br-fourier} H. Brezis:
Equations et in\'equations non lin\'eaires dans les espaces vectoriels en
dualit\'e;
Ann. Inst. Fourier (Grenoble), 18 (1968), 115--175.


\bibitem{croce_rendiconti}
G. Croce: {The regularizing effects of some lower order terms in an
elliptic equation with degenerate coercivity};
{  Rendiconti di Matematica} 27 (2007), 299--314.


\bibitem{ll} J. Leray, J.L. Lions: Quelques r\'{e}sultats de Vi\v{s}ik sur
les
probl\`{e}mes elliptiques semi-lin\'{e}aires par les m\'{e}thodes de Minty et
Browder; {  Bull. Soc. Math. France}, {  93} (1965), 97--107.

 
\bibitem{Po} A. Porretta:  {  Uniqueness and homogenization for a class  ofnon coercive operators in divergence form}; Atti Sem. Mat. Fis. Univ. Modena 46 (1998),   915--936. 


\end{thebibliography}
\end{document}